\documentstyle[12pt]{article}
\newcommand{\nnb}{\nonumber}
\newcommand{\la}{\label}

\newcommand{\be}{\begin{equation}}
\newcommand{\e}{\end{equation}}
\newcommand{\ba}{\begin{eqnarray}}
\newcommand{\ea}{\end{eqnarray}}
\newcommand{\n}{\nabla}

\newcommand{\va}{\varepsilon}
\newcommand{\w}{\infty}

\newcommand{\f}{\frac}
\newcommand{\om}{\Omega}
\newcommand{\p}{\partial}

\begin{document}

\title {Properties of Translating Solutions   to Mean Curvature
 Flow  \thanks{Supported by Natural Science Foundation of China (No. 10631020  and 10871061)
 and the Grant for Ph.D Program
 of Ministry of Education of China. }}

\author{Changfeng Gui\\
{\small Department of Mathematics}\\
 {\small University of Connecticut, Storrs, CT
06269-3009,USA} \\
{\small(e-mail: gui@math.uconn.edu)}\\
Huaiyu Jian and Hongjie Ju\\
  {\small Department of  Mathematics}\\
 {\small Tsinghua University, Beijing 100084, P.R.China}\\
{\small (e-mail: hjian@math.tsinghua.edu.cn)}}

\date {}

\maketitle

{\bf Abstract.}  In this paper, we study  the convexity, interior
gradient estimate, Liouville type theorem and asymptotic behavior at
infinity of translating solutions to mean curvature flow as well as
the nonlinear flow  by powers of the mean curvature.
\vskip0.2cm



\section*{1. Introduction  and Main Results}
\setcounter{section}{1} \setcounter{equation}{0}

In this paper, we study the convexity, interior gradient estimate,
Liouville type theorem and asymptotic behavior at infinity of the
solutions to    equation
\begin{eqnarray} \la {1.1}
a_{ij}u_{ij}:=(\delta_{ij}-\frac{u_iu_j}{1+|\nabla u|^2})u_{ij}
=(\frac{1}{\sqrt{1+|\nabla u|^2}})^{\alpha-1}, \forall x \in R^n
\end{eqnarray}
where constants $\alpha>0$ and we have used the notation
$u_i=\frac{\partial u}{\partial x_i}$ and the convention for
summing.  It is called the {\sl translating soliton equation} of the
nonlinear evolution flow  of hypersurfaces by powers
($\frac{1}{\alpha}$) of the mean curvature. This nonlinear flow was
studied in [10, 11] and  has important applications in minimal
surfaces [2] and isoperimetric inequalities [11]. When $\alpha =1,$
(1.1) is reduced to \be \la {1.2}div
(\frac{Du}{\sqrt{1+|Du|^2}})=\frac{1}{\sqrt{1+|Du|^2}}\ \ in \ \
R^n, \e  which   plays a key role in classifying the type
II-singularity of mean curvature flows [4,5,15]. Scaling the space
and time variables in a proper way near type II-singularity points
on the surfaces evolved by mean curvature vector with a mean convex
initial surface, Huisken-Sinestrari [4,5] and White [15] proved that
the limit flow can be represented as $M_t=\{ (x, u(x)+t)\in R^{n+1}:
x\in R^n, t\in R\}$ and  is also a solution (called a translating
solutions or solitons) to mean curvature flow. Equivalently,  $u$ is
a solution to equation (1.2). Therefore, the classification of type
II-singularity of mean curvature flow is reduced to the
classification of solutions of equation (1.2). There was a well
known conjecture among the geometric flows researchers which asserts
that any complete strictly convex solution of (1.2) is radially
symmetric [15].  A few years ago, Wang [13] proved this conjecture
for $n=2$ and found a non-radially symmetric solution of (1.2) for
$n>2.$     Sheng and Wang [12] used a direct argument   to study the
Singularity profile in mean curvature flow.

One natural question is that how are about the asymptotic behavior
at infinity of the solutions to (1.2), or more generally to equation
(1.1). Obviously, the first step is to make clear the asymptotic
behavior of radially symmetric solutions of (1.1). This leads us to
prove the following theorem 1.1
 in Section 2.

{\bf Theorem 1.1}  { \sl Equation (1.1) has a unique solution of the
form $u(x)=r(|x|)$ up to a translation in $R^{n+1}.$ Moreover, the
function $r\in C^2[0, \infty )$ satisfies that \be \la {1.3}
r''(t)>0,\ \ and \ \ \frac{t}{n}<r'(t)\bigl(1+ (r'(t))^2
\bigr)^{\frac{\alpha-1}{2}}<\f {t}{n-1}\e
 for all $t>0,$
 and
\be \la {1.4}
 r(t)=\f {t^2}{2(n-1)}-\ln t +C_1 - \f
{(n-1)(n-4)}{2}t^{-2}+o(t^{-2}) \quad \hbox{if } \, \alpha =1 \e \be
\la {1.5}
 r(t)=\frac{\alpha}{\alpha+1}\bigl( \f {1}{n-1}
\bigr)^{1/\alpha} t^{1+\frac{1}{\alpha}} - C(\alpha, n)
t^{1-\frac{1}{\alpha}}  +o( t^{1-\frac{1}{\alpha}} ), \quad \hbox
{if } \, \alpha  \not =1 \e
as $t\to \w ,$ where $C_1$ is a constant
depending on $r(1)$ and
$$
C(\alpha, n)=\frac{1}{\alpha-1}  (n-1)^{1/\alpha} \bigl(
\frac{1}{\alpha (n-1)}+ \frac{\alpha -1}{2 }\bigr).
$$
 }

We should mention that when $\alpha=1$,  (1.3) was proved in [8] and
a asymptotic result similar to (1.4) was proven in [1]. Our  method
for general $\alpha>0$ is different and yield more properties of the
solutions. See Section 2 for details.

Another natural question, formulated explicitly as an open problem
in [13], is whether any solution of (1.2) is strictly convex. We
will prove the following theorem 1.2 which is related to  this
question in Section 3.  \vskip0.2cm

{\bf Theorem 1.2}  {\sl Let $u\in C^2(R^n)$ be a convex solution of
equation (1.1). If $u$ is
 strictly convex in some nonempty set, then $u$ is is strictly convex in $R^n.$
 Particularly, $u$ is is strictly convex in $R^n$ if $u(x)\to \infty $ as $|x|\to \infty ;$
 and when $\alpha =1$ (i.e., $ u$ is a solution to (1.2)) and $u(x)\to \infty $ as $|x|\to \infty ,$ then after a rotation of coordinate system,
  $\lim_{h\to \infty}h^{-2}u(hx)\to \sum_{i=1}^kx_i^2$ in $R^n$ for some $k\geq 2,$ and $u$
  is radially symmetric if $n=2.$}

\vskip0.2cm
  This theorem generalizes the main results in [6] which asserts
  that the Hessian $(D^2 u(x))$ has constant rank for all $x\in R^n$
  if $(D^2 u(x))$ is positive semi-definite and $ \Delta u=f(u,
  \nabla u)$ in $R^n,$ where $f\in C^{2, \alpha}$ is strictly positive
  and convex in $u.$ In the case of Minkowski space [9], similar convexity result
  was proved by the second author in [7].

 Theorem 1.1 tells us that the radially symmetric solution of (1.1) is of the $1+\frac{1}{\alpha}$ order growth at infinity.
 This order tends to 1 as $\alpha $ goes to $\w $. Motivated by this, we   have the  following Liouville
result. \vskip0.2cm

 {\bf Theorem 1.3} {\sl There is no nonnegative solution
$u\in C^3(R^n)$ of (1.1) such that
\begin{eqnarray}
\lim_{|x|\rightarrow\infty}\frac{|u(x)|}{|x|}=0.\nonumber
\end{eqnarray}}
\vskip0.2cm

 We will prove this theorem in section 5. For this purpose, in section 4 we will
 want to use
 the gradient estimate techniques by Xu-Jia Wang in [14] and the methods in [3] to
 obtain the following interior gradient result for  equation (1.1).
\vskip0.2cm

 {\bf Theorem 1.4} {\sl Suppose $u\in C^3(B_r(0))$ is a
nonnegative solution of (1.1), then
$$|\nabla u(0)|\leq\exp\{C_1+C_2\frac{M^2}{r^2}\},$$
where $M=\sup_{x\in B_r(0)}u(x),$\ $C_i$(i=1,2) are constants
depending only on $n$ and $\alpha.$ }

      \section*{2. Asymptotic Behavior - Proof of Theorem 1.1 }
\setcounter{section}{2} \setcounter{equation}{0}  We will use a few
lemmas to prove theorem 1.1. The main difficulty is to prove
asymptotic expansion (1.4) and (1.5).

{\bf Lemma 2.1}  { \sl Suppose that $u(x)=r(|x-x_0|)+u(x_0)$ for
some $x_0\in R^n$ and for all $x\in R^n.$  Then $u\in C^2(R^n)$ is a
solution of (1.1) if and only if $r\in C^2(0, \w )$ satisfies \be
\la {2.1} \f {r''}{1+(r')^2}+\f {n-1}{t}r'= \bigl(1+(r')^2
\bigr)^{\frac{1-\alpha}{2}}, \forall t\in (0, \w )\e and \be \la
{2.2} r(0)=r'(0)=0.\e}

{\bf Proof.}  We may assume $x_0=0.$ Let $e_i$ be the unit vector in
positive $x_i$-axis. Since $r(t)=u( t e_i)-u(0)=u(-t e_i)-u(0)$ for
all $t\geq 0 ,$ then $r\in C^2[0, \w )$  and it satisfies  (2.2) if
and only if $u\in C^2(R^n). $
 Writing (1.1) in $r,$ we see that (1.1) is equivalent to (2.1).

 {\bf Lemma 2.2}  {\sl  If $y\in C^1(-\w , \w )$ satisfies
 \be \la {2.3}
 y' +[ (n-1)y (1+y^2)^{\frac{\alpha-1}{2}}-e^s](1+y^2)^{\frac{3-\alpha}{2}}=0, \forall s\in (-\w , \w )
 \e
 and
 \be \la {2.4}\lim_{s \to -\w } \f {y(s)(1+y^2)^{\frac{\alpha-1}{2}}}{e^s}=\f {1}{n}, \e
 then $r(t)=\int_0^t y(\ln s) ds \in C^2[0, \w )$ satisfies (2.1) and (2.2).}

 {\bf Proof.} (2.1)  can be verified directly  by (2.3). Note that (2.4) implies
 $$r'(0):=\lim_{t\to 0^+}r'(t)=\lim_{s \to -\w } y(s) =0=r(0)$$
 and
 $$r''(0):=\lim_{t\to 0^+}r''(t)=\lim_{s \to -\w } \f {y'(s) }{e^s}=\f {1}{n}$$
 by equation (2.3) as well as equation (2.1).

  {\bf Lemma 2.3}  {\sl There exists a  $y\in C^1(-\w , \w )$ solving (2.3) and (2.4) such that
  \be \la {2.5}
  y'(s)>0 \ \ and \ \ \f {e^s}{n}<y(s)(1+y^2)^{\frac{\alpha -1}{2}}<\f {e^s}{n-1}, \forall s\in (-\w , \w ).\e
  Moreover, the function $z(s):= (n-1)e^{-s}y(s)(1+y^2)^{\frac{\alpha-1}{2}}-1$ satisfies
  \be \la {2.6}
  \lim_{s\to -\w} z(s)=-\f {1}{n}, \ \ \lim_{s\to \w} z(s)=0\ \ and \ \ z'(s)>0, \forall s\in (-\w , \w ).
\e} {\bf Proof.} By local existence we have a function $r  \in
C^2[0, \va )$ which solves (2.1) and (2.2) in $(0, \va )$ for some
$\va >0$ (see [8]). Then $y(s)=r'(t) , t=e^s $ satisfies (2.3) in
$(-\w , \ln \va )$ and \be \la {2.7} \lim _{s\to -\w } y(s)=0.\e Let
$(-\w , T)$ be the maximal interval for which $y$ solves (2.3).
First we prove $T=\w .$  In fact, we shall show  that
 \be \la {2.8}
  y'(s)\geq 0 \ \ and \ \  y(s)(1+y^2)^{\frac{\alpha-1}{2}}\leq \f {e^s}{n-1}, \forall s\in (-\w , T ).\e

For convenience,  we may define
$$
g(y)=y(1+y^2)^{\frac{\alpha-1}{2}}.
$$
Note that $\alpha >0$ and
$$
g'(y)=y'(s)(1+y^2)^{\frac{\alpha-3}{2}} (1+\alpha y^2)>0.
$$

 Suppose that  there exists $s_0\in (-\w , T)$ such that
 $ y'(s_0)<0 $ and then  $g(y(s_0))>\f {e^{s_0}}{n-1}$ by (2.3).
We  claim that
   \be \la {2.9} y'(s)<0 ,\ \ and\ \  g(y(s))>\f {e^s}{n-1},  \forall s\in (-\w , s_0 ).\e
   Otherwise, there exists a $s_1<s_0$ such that
   $$ y'(s_1)=0 ,\ \   g(y(s_1))=\f {e^{s_1}}{n-1}\ \ and \ \  y'(s)<0 ,\ \ \forall s\in (s_1 , s_0 ).$$
   Then  $ \f {e^{s_1}}{n-1} =g(y(s_1))>g(y(s_0))>\f {e^{s_0}}{n-1}.$
   This contradiction implies that (2.9) holds, and hence $\lim _{s\to -\w }g(y(s))>\f {e^{s_0}}{n-1}, $
   which  contradicts (2.7). Therefore we have proven (2.8), and then $T=\w
   .$ follows by the standard existence theory of ordinary
   differential equations.

   Next, we  shall prove (2.6). Since $y$ satisfies (2.3) and (2.8) in $(-\w , \w),$ the
    function $z(s):= (n-1)e^{-s}g(y(s))-1$ satisfies
\be \la {2.10} z' +nz +1+  \alpha (n-1)z y^2=0 \ \ and \ \ z\leq 0,
\forall s\in (-\w , \w ).\e Thus, $(ze^{ns})'\geq -e^{ns}$ for all
$s\in (-\w , \w ).$   Integrating this inequality over $(-\w , s)$
and then using (2.7) we obtain \be \la {2.11} z(s)\geq -\f {1}{n},
\forall s\in (-\w , \w ).\e Now we make two observations. The first
one is that $z$ has no local maxima. Indeed,  if $z'(s_0)=0$  for
some $s_0$, we can obtain
 $z''(s_0)> 0$  since
 $z$ satisfies \be \la {2.12} z''+n z' + \alpha (n-1) z' y^2+2
\alpha (n-1) zy y'=0 , \forall s\in (-\w , \w )\e.

 The second observation is  that $z'(s)>0$ for all $s\in (-\w , \w).$
 Otherwise,  if $z'(s_0)\leq 0$ for some $s_0.$
Then  $z'(s)\leq 0$ for all $s\leq s_0, $ since $z$ has no local
maxima. This, together with (2.11), implies that $\lim_{k\to -\w }
z'(s_k)= 0$  for some sequence $s_k\to -\w .$ Thus, $\lim_{k\to -\w
} z(s_k)=-\f {1}{n}$ by (2.10). Since (2.10) and (2.11) imply
$z'(s)<0$ in $(-\infty , s_0),$  we see that $z(s)<-\frac{1}{n}$ for
$s\in (-\infty, s_0),$  which  contradicts (2.11).

Now that $z'(s)>0$ for all $s\in (-\w , \w),$ by (2.10) and (2.11)
we have $\lim_{s\to -\w } z(s)=-\f {1}{n}$ and $\lim_{s\to  \w }
z(s)=0 .$ This proves (2.6),
 which, together with equation (2.3), implies (2.5). The proof of Lemma 2.3 has been completed.

{\bf Lemma 2.4} {\sl Let $y$ be a function as in Lemma 2.3. Then
$r(t)=\int_1^t y(\ln s) ds $ satisfies (1.4) and (1.5) .}

{\bf Proof.} It follows from (2.6) that  $\lim_{k\to \w } z'(s_k)=
0$ for some sequence $s_k\to \w .$ By this, we claim that \be \la
{2.13}  \lim_{s\to \w } z'(s)= 0. \e   Assume that this is not true,
then there exists
 a sequence  of local maxima $\theta_k\to \w $  of $z'(s)$.  Note that we have $z''(\theta _k)=0.$
 From (2.6), we can derive
 $$g(y(s))=\f{e^s}{n-1}  (1+o(1)) \ \ and \ \
 y(s)=(\frac{1}{n-1})^{1/\alpha} e^{s/\alpha}  (1+o(1))
 $$
   as $ s \to \infty .$
 By (2.10) and (2.6), we have
 $\lim_{s \to \w } z'(s)e^{-2s/\alpha}= 0.$
 By the definition of $z$, we have
 $$
 2\alpha (n-1) z y y'= \f{ 2\alpha (n-1)z e^s [ z' +
 (z+1) ] y}{ (1+y^2)^{\f{\alpha-3}{2}} (1+\alpha y^2)} =2(n-1)^{2-2/\alpha} z  [ z' +
 (z+1) ] e^{2s/\alpha} (1+o(1)).
 $$
 This, together with (2.12) and (2.6) again, implies
 $\lim_{k\to \w } z'(\theta _k)=0$. The claim is hence  proven.

 Then it follows from (2.13), (2.10) and (2.6) that
 \be
 z(s)=-\f{1}{\alpha(n-1)}y^{-2}(1+o(1))=-\f{1}{\alpha} (n-1)^{\f{2}{\alpha}-1} e^{-\f{2s}{\alpha}} (1+o(1)) ,
   \ \ as \ \ s\to \w .\e
   and hence
   \be
   g(y)=\frac{1}{n-1} e^s \bigl(1-\f{1}{\alpha} (n-1)^{\f{2}{\alpha}-1} e^{-\f{2s}{\alpha}} \bigr)(1+o(1)),
   \ \ as \ \ s\to \w .
   \e
Therefore, by straightforward computation we obtain
   \be \la {2.14}
   y= e^{\f {s}{\alpha}}
   \bigl( \bigl(\frac{1}{n-1}\bigr)^{1/\alpha}-B(\alpha, n) e^{-\f{2s}{\alpha}} \bigr) (1+o(1)),
   \ \ as \ \ s\to \w ,
   \e
where $$ B(\alpha, n)=(n-1)^{1/\alpha} \bigl( \frac{1}{\alpha^2
(n-1)}+ \frac{\alpha -1}{2 \alpha}\bigr).
$$  This implies (1.5).

In particular,  when $\alpha =1$, we have
$$
y=\frac{1}{n-1}e^s -e^{-s} +o(e^s).
$$

In theory, we can repeat the  above procedure  to obtain higher
order expansion.  Take the simple case  $\alpha =1$ as example, we
can let $w(s)=-\f {e^{2s}}{n-1}z(s)-1.$ Then
 \be \la {2.15}
 \lim_{s\to \w }w(s)=0,\e
 and equation (2.10) read as
 \be \la {2.16}
 w' +(n-2)w+n-2 +\f {e^{2s}}{n-1}w-2(1+w)^2+(n-1)e^{-2s}(1+w)^3=0.\e
 Thus,
 \be \la {2.17}
 w'' + w'[(n-2) +\f {e^{2s}}{n-1} -4(1+w)+\f{3(n-1)}{e^{2s}}(1+w)^2]+\f {2e^{2s}}{n-1}w -
 \f{2(n-1)}{e^{2s}}(1+w)^3 =0.\e
 Since (2.15) means $\lim_{k\to \w }w'(s_k)=0$ for some sequence $s_k\to \w, $ we claim that
 \be \la {2.18}
 \lim_{s\to \w }w'(s)=0.\e
 If not, then there exists a $\delta >0 $ and a sequence $\theta_k\to \w$
  such that $w''(\theta _k)=0$ and
$w'(\theta _k)>\delta$  (or $w'(\theta _k)<-\delta$.)  But it
follows from (2.15) and (2.16) that \be \la {2.19} \lim_{k\to \w }
[w'(\theta _k)+\f{1}{n-1}e^{2\theta_k}w(\theta _k)]= 4-n.\e Taking
 $s=\theta_k$ in (2.17) and using (2.15) and (2.19),
 we have
$$\lim_{k\to \w }[ \f {n-4 }{n-1}w(\theta _k)-\f {4(w^2(\theta _k)+w(\theta _k))}{n-1}+\f {w'(\theta _k)}{n-1}]
 e^{2\theta _k}
=(4-n)(n-2)-4(4-n),$$ which implies $\lim_{k\to \w }
w'(\theta_k)=0,$   a contradiction.

   Thus, using (2.15) and (2.18) we can rewrite (2.16) as
  $$ w=-(n-1)(n-4)e^{-2s}(1+o(1))\ \ as \ \ s\to \w.$$
  By ( 2.14), we have
  $$z(s)=-(n-1)e^{-2s}+(n-1)^2(n-4)e^{-4s}+o(e^{-4s})\ \ as \ \ s\to \w,$$
  which implies
  $$y(s)= \f {e^s}{n-1}-e^{-s}+(n-1)(n-4)e^{-3s}+o(e^{-3s}) \ \ as \ \ s\to \w . $$
  Therefore,  $r(t)=\int_1^t y(\ln s) ds $ satisfies (1.4).

  {\bf Proof of Theorem 1.1:} the existence follows from Lemmas 2.1-2.4; while the
  uniqueness from the well-known comparison principle.

 \section*{3. Convexity - Proof  of Theorem 1.2   }
\setcounter{section}{3} \setcounter{equation}{0}

  {\bf Lemma 3.1} {\sl  Let $u\in C^2(R^n).$  Suppose that there is a constant $c$
   such that the set $\om _c=\{x\in R^n : u(x)<c\}$ is   nonempty and
   bounded. If $u=c$ on $\p \om _c ,$
   then the Hessian matrix $(u_{ij}(x_0))>0$  (positive definite) for some $x_0\in \om _c .$}

  {\bf Proof.}  By the assumption  we see that $u=c$ on $\p \om _c$ and
  $$u(x_1)=\min_{\overline{\om _c} }u(x)<c$$
  for some $x_1\in \om _c.$   This implies that  the function
  $$U(x)=u(x)-\f {c-u(x_1)}{2(diam (\om _c))^2} |x-x_1|^2$$
  must attain interior minimum in $\om _c.$ Consequently, $(u_{ij}(x_0))>0$ for some $x_0\in \om _c ,$
  which implies the desired result.

   {\bf Lemma 3.2} {\sl  Let $u\in C^2(R^n)$ be a convex solution of equation (1.4).  If
    the set   $\om _0 =\{x\in R^n : (u_{ij}(x))>0\}$  is nonempty,  then $  \om _0=R^n .$}

{\bf Proof.} We  follow the arguments of theorem 1.3 by the second
author in [7]. Suppose the contrary that there exists a $x_1\in R^n
\backslash \om _0.$ We will derive a contradiction. We may assume
$\om _0$ is nonempty and connected. (Otherwise, we replace it by one
of its connected components ). Then there exists a short segment
 $l\subset \om _0$ such that $\bar {l}\cap \p \om _0=\{x_1\} .$   Take $x_2\in l$  and $\va >0$
 such that $\overline{B_{\va } (x_2)}\subset \om _0.$  Translating the ball  $B_{\va } (x_2)$ along the
 line $l$  toward $x_1$ we come to a point $\bar{x}$ where the ball and $\p \om _0$
 are touched at the first time. It follows that
 \be \la {3.1}\bar{x}\in R^n \backslash\om _0, \ \ \   B_{\va } (x_0)\subset \om _0 \ \ \
 and \ \ \ \overline{ B_{\va } (x_0)}
 \cap \p \om _0=\{\bar{x}\}\e
 for some $x_0\in \om _0.$ Moreover,  the minimum eigenvalue $\lambda (x)$ of
 the Hessian $(u_{ij}(x))$ satisfies $\lambda (\bar{x})=0.$
  By a coordinate translation and rotation we may
 arrange that
 \be \la {3.2} \bar{x}=0, \ \ u(0)=0, \ \ \n u(0)=0 \ \ and \ \ u_{11}(0)=\lambda (0)=0.\e
 Thus, the origin $0\in  \p B_{\va } (x_0)  $ and
 \be \la {3.3} (u_{ij}(x))>0\ \ in \ \  B_{\va } (x_0).\e
  Rewrite  equation (1.4)   as
 \be \la {3.4} \Delta u =  A(|\n u|^2)u_i u_j u_{i j}+B(|\n u |^2)\ \ in \ \ R^n ,\e
 where $A(t)=\f {1}{t+1}$ and $B(t)=(1+t)^{\f {1-\alpha}{2}}$ are both analytic for $t>-1.$
 Differentiating (3.4) twice with respect to $\f {\p}{\p x_1},$  we have
 \begin{eqnarray}\la {3.5}
\Delta u_{11}&=& 4[A'' u_i u_j u_{i j}+B'' ]u_l u_{l 1}u_m u_{m1}+2[A'  u_{i} u_j u_{i j}+B']u_{m 1} u_{m 1}\nnb \\
 & +& 2[A'  u_{i} u_j u_{i j}+B']u_{m } u_{m 11}+ 8A' u_{m } u_{m 1}u_{i1} u_j u_{i j}\nnb \\
&+& 4A' u_{m } u_{m 1}u_{i} u_j u_{i j1}  +2A u_{i11} u_j u_{i j}\nnb \\
&+& 2A  u_{i1} u_{j1} u_{i j} +4A u_{i1} u_{j} u_{i j1}\nnb\\
&+& A u_{i} u_{j} u_{i j11} \ \  in \ \  R^n .
\end{eqnarray}
 Since $u$ is analytic in $R^n ,$ we expand $u_{11}$ at $x=0$ as a power series to obtain
 $u_{11}(x)=P_k(x)+R(x)$ for all $x\in B_{\delta}(0)$ for some $\delta>0$ such that $\overline {B_{\va }(x_0)}\subset
 B_{\delta}(0) $   (one can choose a smaller $\va $ in advance if necessary),
 where $ P_k(x)$ is the lowest order term, which, by (3.2) and (3.3),  is a nonzero homogeneous
 polynomial of degree $k$, and $R(x)$ is the rest. Note that $k\geq 2$ by the convexity of $u.$
  It follows from (3.3) that  $u_{i i}u_{11}-(u_{i1})^2>0$ in
 $B_{\va }(x_0).$
 Summing over $i$  we have
 \be \la {3.6} \Delta u u_{11}>\sum_{j=1}^n u_{j1}^2\geq  u_{i1}^2\e
 for each $i=1, 2, \cdots , n.$

 We claim that each $u_{i 1}$ is of order at least $\f {k}{2}.$  Otherwise, we expand $u_{i1}$ at $x=0$ as a
  power series so that the lowest order term $h(x)$ must be a a nonzero homogeneous polynomial.
 Choose
 $$a=(a_1, a_2, \cdots , a_n)\in B_{\va }(x_0)\backslash \{x\in B_{\va }(x_0): h(x)=0\}$$
 so that the  segment
 $$L=\{ ta:  t\in (0, 1)\} \subset
  B_{\va }(x_0).$$    Now restricting (3.6) on $L,$ multiplying  the both sides by $t^{-k}$ and then
  letting $t\to 0^+,$  we see the limit of the left-hand side of (3.6) is a nonzero constant multiplied
  by $\Delta u(0)$ which equals to  $1$ by (3.4), but the limit of the right-hand side is positive infinite.
  This is a contradiction.

  Therefore,  each $u_{i 1}$ is of order at least $\f {k}{2}.$
  Hence $u_{ij1},$
  $u_{11i}$ and $ u_{11ij}$ are of order at least $\f {k}{2}-1 ,$ $k-1$ and $k-2$ respectively.  Also note
   that each $u_{i}$ is of order at least 1 by (3.2). With these facts one can check that the right-hand
   side of equation (3.5) is of order at least of $k ;$ while the left-hand side, $\Delta u_{11},$ is
   either of order $k-2,$ or $\Delta P_k=0$ for all $x\in  B_{\va }(x_0).$  Since the first case is
   impossible by comparing  the orders of the two sides, we obtain that $P_k$ is a harmonic polynomial
    in $  B_{\va }(x_0).$

  We claim that  $P_k\geq 0$ for all $x\in  B_{\va }(x_0).$    Otherwise, there exists $a=(a_1, a_2, \cdots ,
   a_n)\in B_{\va }(x_0)$ such that $P_k(a)<0.$ Then
  $$\frac{u_{11}(ta)}{t^k}=P_k(a)+\frac{R(ta)}{t^k}, \ \ \forall t\in (0, 1),$$
  which implies $\lim_{t\to 0^+}\frac{u_{11}(ta)}{t^k}=P_k(a)<0$ contradicting the fact that
  $u_{11}>0$ in $B_{\va }(x_0) $ (see (3.3)).

  Now we use the strong maximal principle to see that $P_k>0$ for all $x\in  B_{\va }(x_0).$   But $P_k(0)=0,$
 and  it follows from Hopf's lemma that $\f {\p P_k}{\p \nu}(0)<0,$ where $\nu $ is the unit outward normal
 to the sphere $\p B_{\va }(x_0).$ This means that the degree of $P_k$ is only one, contradicting the
 fact $k\geq 2.$
This contradiction proves  the lemma.

{\bf Proof of Theorem 1.2:}  it is direct from Lemmas 3.1 and 3.2.
Note that  if $u(x)\to \infty $ as $|x|\to \infty ,$ then $u$ is
strictly convex by Lemmas 3.1 and 3.2. Therefore, when $\alpha =1$,
 we use the results in [13] to know that  after a rotation of coordinate system, $\lim_{h\to
\infty}h^{-2}u(hx)\to \sum_{i=1}^kx_i^2$ in $R^n$ for some $k\geq 2$
by [13; Theorem 1.3] and $u$ is radially symmetric if $n=2$ by [13;
Theorem 1.1].

 \section*{4. Interior Gradient Estimate - Proof  of Theorem 1.4    }
\setcounter{section}{4} \setcounter{equation}{0}

     Let $G(x,\xi)=g(x)\varphi(u)\log u_{\xi}(x),$
where $u$ satisfies the hypothesis of theorem 1.4,
$g(x)=1-\frac{|x|^2}{r^2},\,\varphi(u)=1+\frac{u}{M},\,M=\sup_{x\in
B_{r(0)}}u(x).$  Suppose that $\sup\{G(x,\xi),\,x\in B_r(0),\xi\in
S^{n-1}\}$ is attained at point $x_0$ and in the direction $e_1$.
Then at $x_0,$\,$u_i(x_0)=0$ for $i\geq 2$ since directive
derivatives attain the maximum along the gradient direction, and so
 $a_{11}=\frac{1}{1+u_1^2},\,a_{ii}=1$ for $i\geq2,$ $a_{ij}=0$ for
 $i\neq j.$ As the arguments from (1.2) to (1.4) in [14],  at $x_0$ we have
\begin{eqnarray}\label {4.1}
0=(\log
G)_i=\frac{g_i}{g}+\frac{\varphi^{\prime}}{\varphi}u_i+\frac{u_{1i}}{u_1\log
u_1}
\end{eqnarray}
and
\begin{eqnarray}
0\geq a_{ii}(\log G)_{ii} &\geq & \frac{f_1}{u_1\log
u_1}+\frac{\varphi^{\prime}}{\varphi}f+\frac{u_{11}^2}{2(1+u_1^2)^2\log
u_1}-\frac{2n}{gr^2}-\frac{4}{Mr},\nonumber
\end{eqnarray}
 Where
$$f(x)\equiv (\sqrt{1+|\nabla u(x)|^2})^{1-\alpha}, \ \ f_1(x)\equiv f_{x_1}(x).$$
In particular,
 $$ f(x_0)=(\sqrt{1+u_1^2(x_0)})^{1-\alpha},
\ \ f_1(x_0)= (1-k)(1+u_1(x_0)^2)^{\frac{-\alpha -1}{2}}u_1(x_0)u_{11}(x_0).$$\\
 Suppose $G(x_0, e_1)$ is large enough so that $log u_1>1$ and
$|\frac{g^{\prime}}{g}|\leq\frac{\varphi^{\prime}}{2\varphi}u_1$ at
$x_0,$ then by (4.1) we can obtain $u_{11}\leq
-\frac{\varphi^{\prime}}{2\varphi}u_1^2\log u_1<0.$ Therefore in the
case of $\alpha \geq1$,
\begin{eqnarray}
\frac{f_1}{u_1\log u_1}+\frac{\varphi^{\prime}}{\varphi}f =
\frac{(1-\alpha)u_{11}}{\log
u_1}(1+u_1^2)^{\frac{-\alpha-1}{2}}+\frac{\varphi^{\prime}}{\varphi}(1+u_1^2)^{\frac{1-\alpha}{2}}\geq0\nonumber
\end{eqnarray}
and
\begin{eqnarray}\label{4.2}
\frac{u_{11}^2}{(1+u_1^2)^2}&=&\frac{u_1^2\log^2u_1}{(1+u_1^2)^2}(\frac{g^{\prime}}{g}+\frac{\varphi^{\prime}}{\varphi}u_1)^2\nonumber\\
&\geq&\frac{u_1^4\log^2u_1}{(1+u_1^2)^2}(\frac{\varphi^{\prime}}{2\varphi})^2\nonumber\\
&\geq&\frac{{\varphi^{\prime}}^2}{8\varphi^2}\log^2u_1\nonumber\\
&\geq&\frac{\log^2u_1}{32M^2}.
\end{eqnarray}\\
 If $0<\alpha<1,$ then there exists a positive integer $m$ such that
$\alpha m>1,$ we may suppose that $G(x_0, e_1)$ is still suitably
large so that $log u_1>1$ and
$|\frac{g^{\prime}}{g}|\leq\frac{\varphi^{\prime}}{m\varphi}u_1$ at
$x_0$,  then by (4.1) we have
\begin{eqnarray}\label{4.3}
\frac{f_1}{u_1\log u_1}+\frac{\varphi^{\prime}}{\varphi}f & =
&\frac{(1-\alpha)u_{11}}{\log
u_1}(1+u_1^2)^{\frac{-\alpha-1}{2}}+\frac{\varphi^{\prime}}{\varphi}(1+u_1^2)^{\frac{1-\alpha}{2}}\nonumber\\
&=&(1+u_1^2)^{\frac{-\alpha-1}{2}}[\frac{(1-\alpha)u_{11}}{\log
u_1}+\frac{\varphi^{\prime}}{\varphi}(1+u_1^2)]\nonumber\\
&=&(1+u_1^2)^{\frac{-\alpha-1}{2}}[-\frac{(1-\alpha)g_1u_1}{g}-\frac{(1-\alpha)\varphi^{\prime}u_1^2}{\varphi}
+\frac{\varphi^{\prime}}{\varphi}(1+u_1^2)]\nonumber\\
&=&(1+u_1^2)^{\frac{-\alpha-1}{2}}[-\frac{(1-\alpha)g_1u_1}{g}+\frac{\varphi^{\prime}}{\varphi}(1+\alpha u_1^2)]\nonumber\\
&\geq&(1+u_1^2)^{\frac{-\alpha-1}{2}}\frac{\varphi^{\prime}}{m\varphi}[(\alpha-1)u_1^2+m+m\alpha
u_1^2]\geq0.
\end{eqnarray}\\
 Obviously, in this case (4.2) becomes
\begin{eqnarray}
\frac{u_{11}^2}{(1+u_1^2)^2} \geq C\frac{\log^2u_1}{M^2},\nonumber
\end{eqnarray}
where the constant $C$ depends only on $\alpha .$ \\
To sum up the above two cases we see that at $x_0$ and for any
$\alpha>0,$
\begin{eqnarray}
0\geq a_{ii}(\log G)_{ii} &\geq & \frac{f_1}{u_1\log
u_1}+\frac{\varphi^{\prime}}{\varphi}f+\frac{u_{11}^2}{2(1+u_1^2)^2\log
u_1}-\frac{2n}{gr^2}-\frac{4}{Mr}\nonumber\\
&\geq&C\frac{\log u_1}{M^2}-\frac{2n}{gr^2}-\frac{4}{Mr}.\nonumber
\end{eqnarray}
Recall  we may have assumed $\log u_1(x_0)>1.$  Then we obtain
\begin{eqnarray}
g(x_0)\log u_1(x_0) \leq C_3\frac{M}{r}+C_4\frac{M^2}{r^2},\nonumber
\end{eqnarray}
where $C_3$ and $C_4$ depend only on $n$ and $\alpha$.\\
 Since
$ 1\leq\varphi\leq 2 ,$ we use   Cauchy inequality to get
\begin{eqnarray}
G(x_0,e_1)&=&g(x_0)\varphi(u(x_0))\log
u_1(x_0)\nonumber\\
&\leq&2C_3\frac{M}{r}+2C_4\frac{M^2}{r^2}\nonumber\\
&\leq &C_1(n,\alpha)+C_2(n,\alpha)\frac{M^2}{r^2}.\nonumber
\end{eqnarray}
Therefore, we have proved that for any $\xi\in S^{n-1},$
\begin{eqnarray}
G(x_0,e_1)\geq G(0,\xi)\geq\log u_{\xi}(0),\nonumber
\end{eqnarray}
 which implies
\begin{eqnarray}
u_{\xi}(0)\leq\exp\{C_1+C_2\frac{M^2}{r^2}\}.\nonumber
\end{eqnarray}
Noting that $\xi $ can be any vector in $S^{n-1},$  we have
$$|\nabla u(0)|\leq\exp\{C_1+C_2\frac{M^2}{r^2}\}.$$
This proves theorem 1.4.

 \section*{5. Liouville Type Theorem - Proof  of Theorem 1.3   }
\setcounter{section}{5} \setcounter{equation}{0}

  Suppose $u\in C^3(R^n)$ is a nonnegative solution
of (1.4),and
\begin{equation}
|u(x)|=o(|x|),\quad as \quad |x|\rightarrow\infty.
\end{equation}
We will prove $u\equiv constant $. Since any constant is not a
solution to (1.4), the theorem will be proved.

 By Theorem 1.4 and
(5.1) we have
\begin{equation}
|\nabla u(x)|\leq C,\quad \forall x\in R^n.
\end{equation}
We claim that
$$\nabla u(x)=0,\quad \forall x\in R^n.$$
Otherwise, there is a $y$ such that $\nabla u(y)\neq 0$ for some
$y.$ We may assume $y=0$ and thus
\begin{equation}
|\nabla u(0)|\geq\delta>0
\end{equation}
for some positive $\delta . $  We will induce a contradiction.

   Let
$$ G(x,\xi)=g(x)\varphi(u)u_{\xi}(x),$$ where
$$g(x)=1-\frac{|x|^2}{r^2},\,\varphi(u)=(1-\frac{u}{M})^{-\beta},\,M=4\sup\{|u(x)|,x\in
B_r(0)\}$$ and $\beta \in(0,1)$ is to  be determined.  Suppose
$\sup\{G(x,\xi),\,x\in B_r(0),\xi\in S^{n-1}\}$ is attained at point
$x_0$ and in the direction $e_1$.\\
By (5.3) we have
\begin{equation}
g(x_0)\geq\delta_1,\ u_1(x_0)\geq\delta_1
\end{equation}
for some $\delta_1>0$ depending only on $\delta$.\\
Then at $x_0$,
\begin{eqnarray}
0=(\log
G)_i=\frac{u_{1i}}{u_1}+\frac{g_i}{g}+\frac{\varphi^{\prime}}{\varphi}u_i
\end{eqnarray}
and
\begin{eqnarray}
(\log
G)_{ij}&=&\frac{u_{1ij}}{u_1}+(\frac{\varphi^{\prime\prime}}{\varphi}-2\frac{{\varphi^{\prime}}^2}{\varphi^2})u_iu_j
+\frac{\varphi^{\prime}}{\varphi}u_{ij} +(\frac{g_{ij}}{g}
-2\frac{g_ig_j}{g^2})\nonumber\\
&-&\frac{\varphi^{\prime}}{g\varphi}(g_iu_j+g_ju_i)\nonumber
\end{eqnarray} where we have used (5.5).
Note that at $x_0$, $u_i=0$ for $i\geq2,$
$a_{11}=\frac{1}{1+u_1^2},\ a_{ii}=1$ for $i\geq2,$ and $a_{ij}=0$
for $i\neq j.$   Therefore at $x_0,$
\begin{eqnarray}
0\geq a_{ii}(\log G)_{ii}&=&\frac{a_{ii}u_{1ii}}{u_1}+
(\frac{\varphi^{\prime\prime}}{\varphi}-2\frac{{\varphi^{\prime}}^2}{\varphi^2})
\frac{u_1^2}{1+u_1^2}+\frac{\varphi^{\prime}}{\varphi}f\nonumber\\
&+&a_{ii}(\frac{g_{ii}}{g}
-2\frac{g_i^2}{g^2})-\frac{2g_1\varphi^{\prime}u_1}{g\varphi(1+u_1^2)}.
\end{eqnarray}
By (5.4) we obtain
\begin{eqnarray}
-\frac{2g_1\varphi^{\prime}u_1}{g\varphi(1+u_1^2)}&=&\frac{4\beta}{M-u}\cdot\frac{x_1}{gr^2}
\cdot\frac{u_1}{1+u_1^2}\nonumber\\
&\geq&\frac{-16\beta}{3M\delta_1r}\nonumber\\
&=&-\frac{C_1}{Mr}
\end{eqnarray}
and
\begin{eqnarray}
a_{ii}(\frac{g_{ii}}{g} -2\frac{g_i^2}{g^2})&=&-\frac{2}{g
r^2}(\frac{1}{1+u_1^2}+n-1)-\frac{8x_1^2}{(1+u_1^2)g^2r^4}
-\sum_{i=2}^n \frac{8x_i^2}{g^2r^4}\nonumber\\
&\geq&-\frac{C_2}{r^2},
\end{eqnarray}
where constants  $C_1,C_2$ depend only on $n,$ $\delta_1$ and $\beta$.\\
Differentiating   equation (1.4) with respect to $x_1$, we have
$$\sum_{i=1}^n u_{1ii}=f_1+(\frac{u_iu_j}{1+u_1^2}u_{ij})_{x_1}$$
which implies $$\frac{u_{111}}{1+u_1 ^2}=f_1-\sum_{i=2}^nu_{1ii}$$
at $x_0$.  Thus, we  obtain that at $x_0,$
\begin{eqnarray}
\frac{a_{ii}u_{1ii}}{u_1}&=&\frac{1}{u_1}[f_1+\frac{2u_1u_{11}^2}{(1+u_1^2)^2}
+\sum_{i\geq2}\frac{2u_1}{1+u_1^2}u_{1i}^2]\nonumber\\
&\geq&\frac{1}{u_1}f_1\nonumber\\
&=&(1-\alpha)(1+u_1^2)^{-\frac{\alpha+1}{2}}u_{11}.\nonumber
\end{eqnarray}
 Observing that (5.5) implies
$u_{11}=-\frac{u_1^2\varphi^{\prime}}{\varphi}-\frac{u_1g_1}{g}$ at
$x_0$, we see that
\begin{eqnarray}
\frac{a_{ii}u_{1ii}}{u_1}+\frac{\varphi^{\prime}}{\varphi}f&\geq&(1-\alpha)(1+u_1^2)^{-\frac{\alpha+1}{2}}u_{11}
+\frac{\varphi^{\prime}}{\varphi}(1+u_1^2)^{\frac{1-\alpha}{2}}\nonumber\\
&=&(1-\alpha)(1+u_1^2)^{-\frac{\alpha+1}{2}}[-\frac{u_1^2\varphi^{\prime}}{\varphi}-\frac{u_1g_1}{g}]
+\frac{\varphi^{\prime}}{\varphi}(1+u_1^2)^{\frac{1-\alpha}{2}}\nonumber\\
&=&\frac{\varphi^{\prime}}{\varphi}(1+u_1^2)^{-\frac{\alpha+1}{2}}[\alpha
u_1^2-u_1^2+1+u_1^2]-
(1-\alpha)(1+u_1^2)^{-\frac{\alpha+1}{2}}\frac{u_1g_1}{g}\nonumber\\
&=&\frac{\alpha}{M-u}(1+u_1^2)^{-\frac{\alpha+1}{2}}(\alpha
u_1^2+1)+2(1-\alpha)(1+u_1^2)^{-\frac{\alpha+1}{2}}\frac{u_1x_1}{gr^2}.\nonumber\\
& \geq & \frac{1}{r}(1+u_1^2)^{-\frac{\alpha+1}{2}}[
\frac{r}{M-u}\alpha (\alpha u_1^2+1)-2|1-\alpha|
\frac{u_1}{\delta_1}]. \nonumber
\end{eqnarray}
  Since $\delta_1<u_1\leq C$ and $\lim_{r\to \infty}
  \frac{r}{M-u}=+\infty$ by (5.1),
  we have for $r$ large enough,
\begin{eqnarray}
\frac{a_{ii}u_{1ii}}{u_1}+\frac{\varphi^{\prime}}{\varphi}f\geq0.
\end{eqnarray}
Inserting (5.7),(5.8) and (5.9)into (5.6), we obtain that at $x_0$
and for large $r$,
\begin{eqnarray}
(\frac{\varphi^{\prime\prime}}{\varphi}-2\frac{{\varphi^{\prime}}^2}{\varphi^2})
\frac{u_1^2}{1+u_1^2}\leq\frac{C_1}{Mr}+\frac{C_2}{r^2}.
\end{eqnarray}
Now chose  $\beta \in (0,1)$   so that
\begin{eqnarray}
\frac{\varphi^{\prime\prime}}{\varphi}-2\frac{{\varphi^{\prime}}^2}{\varphi^2}=\frac{\beta
(1-\beta) }{(M-u)^2}\geq \frac{1}{10 M^2}.
\end{eqnarray}
 Then (5.10) implies
\begin{eqnarray}
\frac{u_1^2}{1+u_1^2}\leq 10 [ C_1\frac{M}{r}+C_2\frac{M^2}{r^2}]
\nonumber
\end{eqnarray} at $x_0$
and for large $r$.  Let $r\rightarrow\infty,$ by (5.1) we obtain
$u_1(x_0)=0,$ contradicting (5.4).  In this way, we have proved
theorem 1.3.

 \newpage
\section* { References}

\begin{enumerate}
\itemsep -2pt

\item[1] Altschuler, Steven; Angenent, Sigurd B.; Giga, Yoshikazu Mean
curvature flow through singularities for surfaces of rotation. J.
Geom. Anal. 5 (1995), no. 3, 293--358.

\item [2] T. H. Colding and W. P. Minicozzi II, Width and mean
curvature flow, preprint, Arxiv: 0705.3827vz, 2007.

 \item [3] D. Gilbarg and N. S. Trudinger, Elliptic partial
differential equations of second order, 2nd Version,
Springer-Verlag, 1983

    \item[4] G. Huisken and C. Sinestrari, Mean curvature flow
singularities for mean convex surfaces, Calc. Var. Partial Differ.
Equations, 8(1999), 1-14.
\item[5]G. Huisken and C. Sinestrari, Convexity estimates for
mean curvature flow and singularities of mean convex surfaces,
Acta Math., 183(1999), 45-70.

\item [6] N. Korevaar and J. Lewis, convex solutions to
       certain equations have constant rank Hessian, arch Rational
       Mech. Anal., 97(1987), 19-32.

 \item[7] H.Y.  Jian,   Translating solitons of mean curvature flow of noncompact
 spacelike hypersurfaces in Minkowski space, J. Differential Equations, 220(2006),
 147-162.

\item [8] H. Y. Jian, Q. H. Liu and X. Q Chen, Convexity and
symmetry of translating solitons in mean curvature flows, Chin. Ann.
Math., 26B (2005), 413-422.

 \item [9] Y. N. Liu and H. Y. Jian, Evolution of spacelike
       hypersurfaces by mean curvature minus external force field in
       minkowski space,
       to appear in  Advanced Nonlinear Studies, Aug 2009.

\item [10] F. Schulze, Evolution of convex hypersurfaces by powers of
the mean curvature, Math Z. 251(2005), 721-733.
\item [11] F. Schulze, Nonlinear Evolution by mean curvature and
isoperimetric inequalities, J. Differential Geom., 79(2008),
197-241.

\item [12] W. M. Sheng and X. J. Wang, Singularity profile in the
mean curvature flow,arXiv: math.DG/0902.2261v2.

 \item [13] X. J. Wang, Translation solutions to the mean curvature flow, arXiv:
 math.DG/0404326v1 (submitted in Ann Math, 2003.)

 \item[14] X. J.Wang,Interior gradient estimates for mean curvature
equations,Math. Z.228 (1998), 73-81.

\item [15] B. White, The nature of singularities in mean
 curvature flow of mean-convex sets, J. Amer. Math. Soc., 16(2003),
 123-138.

 \end{enumerate}

\end{document}